\documentclass[12pt]{article}
\usepackage{amsthm,amssymb,amsmath}
\usepackage{graphicx}

\newcommand{\sect}[1]{\section{#1}\setcounter{equation}{0}}
\newcommand{\subsect}[1]{\subsection{#1}}

\newtheorem{theorem}{Theorem}[section]
\newtheorem{lemma}[theorem]{Lemma}

\newtheorem{proposition}[theorem]{Proposition}

\theoremstyle{definition}

\newtheorem{definition}[theorem]{Definition}

\newtheorem{remark}[theorem]{Remark}

\newcommand{\Log}{{\rm Log}}

\newcommand{\Imag}{{\rm Im}}
\newcommand{\Real}{{\rm Re}}

\begin{document}

\date{}

\title
{Approximation theorems for Banach-valued almost periodic and semi-almost periodic holomorphic functions}


\author{Alexander Brudnyi\thanks{Research supported in part by NSERC.
}\\
Department of Mathematics and Statistics\\
University of Calgary, Calgary\\
Canada\\
\\
Damir Kinzebulatov\\
Department of Mathematics\\
University of Toronto, Toronto\\
Canada}

\maketitle

\begin{abstract}
The paper 
studies semi-almost periodic holomorphic functions on a polydisk which have, in a sense, the weakest possible discontinuities on the boundary torus. The basic result used in the proofs is an
extension of the classical Bohr approximation theorem for almost periodic holomorphic functions on a strip to the case of Banach-valued almost periodic holomorphic functions.
\end{abstract}

\sect{Introduction and Formulation of Main Results}



{\bf 1.1.} Let $B$ be a complex Banach space and $C^B_b(U)$ be  the space of bounded continuous $B$-valued functions $f$ on a closed subset $U\subset\mathbb C$ endowed with $\sup$-norm $\|f\|:=\sup_{t \in U}\|f(t)\|_B$.

\begin{definition}
A function $f\in C^B_b(\mathbb R)$ is said to be \textit{almost periodic} if the family of shifts $\{S_\tau f\}_{\tau \in \mathbb R}$\,, $S_\tau f(t):=f(t+\tau)$, $t\in\mathbb R$, is relatively compact in $C_b^B(\mathbb R)$. 
\end{definition}

By $AP_B(\mathbb R)\subset C^B_b(\mathbb R)$ we denote the Banach space of $B$-valued almost periodic functions on $\mathbb R$ endowed with $\sup$-norm.

The following result due to H.~Bohr, commonly called the \textit{approximation theorem}, is a cornerstone of the theory of Banach-valued almost periodic functions.

\begin{theorem}[H.~Bohr]
\label{thm1}
The vector subspace of $C^B_b(\mathbb R)$ spanned over $B$ by functions $t \mapsto e^{i\lambda t}$, $\lambda \in \mathbb R$, is dense in $AP_B(\mathbb R)$.
\end{theorem}

Similarly one defines $B$-valued holomorphic almost periodic functions on the strip $\Sigma:=\{z \in \mathbb C: 0 \leq \Imag(z) \leq \pi\}$.

\begin{definition}\label{def2}
A function $f\in C^B_b(\Sigma)$ is called \textit{holomorphic almost periodic} on $\Sigma$ if it is holomorphic in the interior of $\Sigma$
and the family of shifts $\{S_x f\}_{x \in \mathbb R}$ is relatively compact in $C_b^B(\Sigma)$. 
\end{definition}

By $APH_B(\Sigma)\subset C^B_b(\Sigma)$ we denote the Banach space of $B$-valued almost periodic holomorphic functions on $\Sigma$ endowed with $\sup$-norm.

Our first result is an analog of Theorem \ref{thm1} for functions from $APH_B(\Sigma)$.
\begin{theorem}
\label{thm3}
The complex vector subspace of $C^B_b(\Sigma)$ spanned over $B$ by functions $z \mapsto e^{i\lambda z}$, $\lambda \in \mathbb R$, is dense in $APH_{B}(\Sigma)$.
\end{theorem}

This result is equivalent to the fact that the Banach space $APH(\Sigma):=APH_{\mathbb C}(\Sigma)$ of holomorphic almost periodic functions on $\Sigma$ has the {\em approximation property} along with the scalar version of Theorem \ref{thm3}, see Grothendieck \cite[Section 5.1]{Gro}. (Similarly Theorem \ref{thm1} can be derived from the fact that $AP(\mathbb R):=AP_{\mathbb R}(\mathbb R)$ has the approximation property.) Let us recall
\begin{definition}\label{def2'}
A Banach space $B$ is said to have the approximation property, if, for every compact set $K\subset B$ and every $\varepsilon>0$, there is an operator $T:B\to B$ of finite rank so that $\|Tx-x\|_B\le\varepsilon$ for every $x\in K$. If such $T$ has norm at most $1$, then $B$ is said to have metric approximation property.
\end{definition}

It is easy to see that every space $B$ with a Schauder basis has the approximation property. However, there are Banach spaces without this property, the first such example was constructed by Enflo \cite{E}. In Proposition \ref{prop3.1} we will show that $APH(\Sigma)$ has metric approximation property.

\noindent {\bf 1.2.}
We apply Theorem~\ref{thm3} to study the algebra of bounded holomorphic functions on the unit polydisk $\mathbb D^n\subset\mathbb C^n$ with semi-almost periodic boundary values. Originally the algebra of semi-almost periodic functions on the real line $\mathbb R$ was introduced and studied by Sarason \cite{Sar} in connection with some important problems of operator theory. By the definition it is a closed subalgebra of $C_b(\mathbb R)$ generated by algebras $AP(\mathbb R)$ and $PC_\infty$ (the algebra of continuous
functions on $\mathbb R$ that have finite limits at $+\infty$ and at $-\infty$).
In \cite{BK} we introduced the algebra $SAP(\partial\mathbb D)$ of semi-almost periodic functions on the unit circle $\partial\mathbb D$ which generalizes the Sarason algebra. Then we studied the algebra $H^\infty(\mathbb D) \cap SAP(\partial\mathbb D)$ of bounded holomorphic functions on the unit disk $\mathbb D$ having boundary values in $SAP(\partial\mathbb D)$ which we called  \textit{the algebra of semi-almost periodic holomorphic functions}.
(Recall that for a function $f \in H^\infty(\mathbb D)$, one defines its boundary values a.e. on $\partial\mathbb D$ by taking limits in non-tangential directions, so that $f|_{\partial\mathbb D} \in L^\infty(\partial\mathbb D)$.) 

Our interest in algebra $H^\infty(\mathbb D) \cap SAP(\partial \mathbb D)$ is motivated by the following problem:
\begin{itemize}
\item[]
{\em Find a closed subalgebra $A\subset H^\infty(\mathbb D)$ whose elements have, in a sense, the weakest possible discontinuities at the boundary $\partial \mathbb D$.}
\end{itemize}
An obvious candidate  
would be the subalgebra of bounded holomorphic functions having discontinuities of at most first kind at the boundary. Unfortunately, as it follows from the Lindel\"{o}f Theorem, see \cite{Gar}, any such function must be continuous on $\bar{\mathbb D}$. Moreover, the same result holds if we consider bounded holomorphic functions with first kind discontinuous on $\partial\mathbb D$ of their real or imaginary parts. 

Another way to measure singularities of a function from $H^{\infty}(\mathbb D)$ is to consider discontinuities of its modulus on $\partial\mathbb D$. Assuming that
$A\subset H^\infty(\mathbb D)$ consists of functions $f$ such that $|f||_{\partial\mathbb D}\in L^{\infty}(\partial\mathbb D)$ has 
first-kind discontinuous only, we obtain that
$A$ contains all inner functions. Then by the Marshall theorem, see \cite{Gar}, $A=H^\infty(\mathbb D)$. Thus to get a nontrivial answer we correct the above problem restricting ourselves to the case of algebras $A\subset H^{\infty}(\mathbb D)$ generating by subgroups $G_A$ of invertible elements of $H^\infty(\mathbb D)$ such that for each $f\in G_A$ the function $|f|$ (or, equivalently, the harmonic function $\ln |f|$) has finitely many discontinuities of at most first kind on $\partial\mathbb D$. The main result in \cite{BK} (Theorem 1.8) establishes a connection between such algebras $A$ and certain subalgebras of $H^\infty(\mathbb D) \cap SAP(\partial \mathbb D)$.

The purpose of the present paper is to define semi-almost periodic holomorphic functions on $\mathbb D^n$ and to
extend the results of \cite{BK} to this class of functions.  

First we recall some definitions from \cite{BK}.

In what follows, we consider $\partial\mathbb D$ with the counterclockwise orientation. For $t_0 \in \mathbb R$ let $$\gamma_{t_0}^k(s):=\{e^{i(t_0+kt)} :\ 0 \leq t<s<2\pi\}, \quad k \in \{-1,1\},$$
be two open arcs having $e^{it_0}$ as the right and the left endpoints, respectively.

\begin{definition}[\cite{BK}]
\label{def1}
A function $f \in L^\infty(\partial\mathbb D)$ is called semi-almost periodic on $\partial\mathbb D$ if for any $t_0 \in [0,2\pi)$ and any $\varepsilon>0$ there exist a number $s=s(t_0,\varepsilon) \in (0,\pi)$ and functions $f_k:\gamma_{t_0}^k(s)\to\mathbb C$, $k\in\{-1,1\}$, such that functions
\begin{equation*}
t\mapsto f_k\bigl(e^{i(t_0+kse^t)}\bigr), \quad -\infty<t<0, \quad k \in \{-1,1\},
\end{equation*}
are restrictions of some almost periodic functions from $AP(\mathbb R)$,
and
\begin{equation*}
\sup_{z \in \gamma_{t_0}^k(s)}|f(z)-f_k(z)|<\varepsilon ,\ \ \ k\in\{-1,1\}.
\end{equation*}
\end{definition}

We denote by $SAP(\partial\mathbb D)$ the Banach algebra of semi-almost periodic functions on $\partial\mathbb D$ endowed with $\sup$-norm. For $S$ a closed subset of $\partial\mathbb D$ we denote by $SAP(S)$ the Banach algebra of semi-almost periodic functions on $\partial\mathbb D$ that are continuous on $\partial\mathbb D\setminus S$. (Note the the Sarason algebra is isomorphic to $SAP(\{z_0\})$ for any $z_0\in\partial\mathbb D$.)

Next, let $A({\mathbb D})$ be the algebra consisting of holomorphic functions in $H^\infty(\mathbb D)$ that are continuous on $\bar{\mathbb D}$. Suppose that $S$ contains at least two points. By $A_S$ we denote the closure in $L^\infty(\partial \mathbb D)$ of the algebra generated by $A(\mathbb D)$ and holomorphic functions of the form $e^{f}$, where $\Real(f)|_{\partial\mathbb D}$ is a finite linear combination (over $\mathbb R$) of characteristic functions of closed arcs in $\partial\mathbb D$ whose endpoints belong to $S$.
If $S$ consists of a single point, then we define $A_S$ to be the closure of the algebra generated by $A(\mathbb D)$ and functions $ge^{\lambda f}$, where $\Real(f)|_{\partial\mathbb D}$ is the characteristic function of a closed arc with an endpoint in $S$ and $g \in A(\mathbb D)$ is a function such that $ge^{f}$ has discontinuity on $S$ only. 

Now, Theorem 1.8 of \cite{BK} describes the structure of the algebra of semi-almost periodic holomorphic functions:
\begin{theorem}[\cite{BK}]
\label{thm4}
\label{oldapprox}
$H^\infty(\mathbb D) \cap SAP(S)=A_S$.
\end{theorem}
\noindent (Here and below we identify the elements of algebra $H^\infty(\mathbb D^n)$ with their boundary values defined on $(\partial\mathbb D)^n$.)
\begin{remark}\label{rem1}
Suppose that $S\subset\partial\mathbb D$ contains at least 2 points. Let $e^{\lambda f}\in A_{S}$, $\lambda\in\mathbb R$, where $\Real(f)$ is the characteristic function of an arc $[x,y]$ with $x,y\in S$. Let
$\phi_{x,y}:\mathbb D\to\mathbb H_{+}$ be the bilinear map onto the upper half-plane that maps $x$ to $0$, the midpoint of the arc $[x,y]$ to $1$ and $y$ to $\infty$. Then
there is a constant $C$ such that
$$
e^{\lambda f(z)}=e^{-\frac{i\lambda}{\pi}{\rm Log}\ \! \phi_{x,y}(z)+\lambda C},\ \ \
z\in\mathbb D,
$$
where ${\rm Log}$ is the principal branch of the logarithmic function. Thus Theorem \ref{thm4} implies that the algebra $H^{\infty}(\mathbb D)\cap SAP(S)$ is the uniform closure of the algebra generated by $A(\mathbb D)$ and the family of functions
$e^{i\lambda ({\rm Log}\circ\phi_{x,y})}$, $\lambda\in\mathbb R$, $x,y\in S$.
\end{remark}

Let us define now multi-dimensional analogues of algebras $SAP(S)$ and $A_S$. Namely, if $S_k \subset \partial\mathbb D$ are closed sets, $1 \leq k \leq n$, and $S=\prod_{k=1}^n S_k \subset (\partial\mathbb D)^n$, then we define
\begin{equation*}
SAP^n(S):=\bigotimes_{k=1}^n SAP(S_k),\ \ \
A^n_S:=\bigotimes_{k=1}^n A_{S_k}.
\end{equation*}

Here $\bigotimes$ stands for completion of symmetric tensor product of the corresponding algebras. In particular, $SAP^n(S)$ and $A_S^n$ are uniform closures of the algebras of complex polynomials in variables $f_{1}(z_1),\dots, f_n(z_n)$ with $f_k(z_k)\in SAP(S_k)$ and $f_k(z_k)\in A_{S_k}$, $1\leq k\leq n$, respectively, where $z=(z_1,\dots, z_n)\in (\partial\mathbb D)^n$.
 
Now the extension of Theorem \ref{thm4} is read as follows.

\begin{theorem}
\label{approx2}
$H^\infty(\mathbb D^n) \cap SAP^n(S)=A^n_S$.
\end{theorem}

\begin{remark}
We apply Theorem \ref{thm3} and some arguments from \cite{BK} to show that the algebra $H^\infty(\mathbb D) \cap SAP(S)$ has  the approximation property. This will imply Theorem \ref{approx2}.
\end{remark}

Let $M_S$ be the maximal ideal space of the algebra $H^\infty(\mathbb D^n) \cap SAP^n(S)$ and $M_{S_k}$ be the maximal ideal space of the algebra $H^\infty(\mathbb D) \cap SAP(S_k)$, $1\leq k\leq n$.
As a corollary of Theorem \ref{approx2} we obtain 

\begin{theorem}\label{max}
$M_S$ is homeomorphic to $M_{S_1}\times\cdots\times M_{S_n}$.
Moreover, the Corona Theorem is valid for $H^\infty(\mathbb D^n) \cap SAP^n(S)$, i.e., $\mathbb D^n$ is dense in $M_S$ in the Gelfand topology.
\end{theorem}

\begin{remark}\label{struct}
The structure of $M_{S_k}$ is described in \cite{BK}, Theorems 1.7 and 1.14. Let us recall this result. In what follows $M(A)$ stands for the maximal ideal space of the Banach algebra $A$.

Since $A(\mathbb D)\hookrightarrow A_{S_k}$, there is a continuous
surjection of the maximal ideal spaces 
$$
a_{S_k}:M_{S_k}\to M(A(\mathbb D))\cong\bar{\mathbb D}.
$$
Recall that the {\em \v{S}ilov boundary} of $A_{S_k}$ is the smallest compact subset $K\subset M_{S_k}$ such that for each $f\in A_{S_k}$
$$
\sup_{z\in M_{S_k}}|f(z)|=\sup_{\xi\in K}|f(\xi)|.
$$
Here we assume that every $f\in A_{S_k}$ is also defined on $M_{S_k}$ where its extension to $M_{S_k}\setminus\mathbb D$ is given by the {\em Gelfand transform}: $f(\xi):=\xi(f)$, $\xi\in M_{S_k}$.
\begin{theorem}[\cite{BK}]
\begin{itemize}
\item[(1)] $a_{S_k}: M_{S_k}\setminus a_{S_k}^{-1}(S_k)\to\bar{\mathbb D}\setminus S_k$ is a homeomorphism.
\item[(2)] The \v{S}ilov boundary $K_{S_k}$ of $A_{S_k}$ is naturally homeomorphic to $M(SAP(S_k))$. Moreover, $K_{S_k}\setminus a_{S_k}^{-1}(S_k)=\partial\mathbb D\setminus S_k$ and $K_{S_k}\cap a_{S_k}^{-1}(z)$, $z\in S_k$, is homeomorphic to the disjoint union $b\mathbb R\sqcup b\mathbb R$ of the Bohr compactifications $b\mathbb R$ of $\mathbb R$.
\item[(3)]
For each $z\in S_k$ preimage $a_{S_k}^{-1}(z)$ is homeomorphic to the maximal ideal space of the algebra $APH(\Sigma)$ of holomorphic almost periodic functions on the strip $\Sigma:=\{z\in\mathbb C\ :\ \Imag(z)\in [0,\pi]\}$.
\end{itemize}
\end{theorem}
\end{remark}

{\em Acknowledgment.} We are grateful to S.\,Favorov for useful discussions and valuable comments improving the presentation.
%

\sect{Auxiliary Results}

{\bf 2.1.} In our proofs we use some results on Bohr's compactifications.

First recall that the algebra of almost periodic functions $AP(\mathbb R):=AP_{\mathbb C}(\mathbb R)$ is naturally isomorphic to the algebra $C(b\mathbb R)$ of complex continuous functions on the Bohr compactification $b\mathbb R$ of $\mathbb R$ (the maximal ideal space of $AP(\mathbb R)$). 
That is, a complex continuous function on $\mathbb R$ is almost periodic if and only if it admits a continuous extension to $b\mathbb R$ by means of the Gelfand transform. (Recall that the maximal ideal space $b\mathbb R$ is defined as the space of continuous non-zero characters $AP(\mathbb R) \mapsto \mathbb C$ endowed with the Gelfand (i.e. weak*) topology. Then $b\mathbb R$ is a compact abelian group, $\mathbb R$ is naturally embedded into $b\mathbb R$ as a dense subset, so that the action of $\mathbb R$ on itself by translations extends uniquely to the continuous action of $\mathbb R$ on $b\mathbb R$.)

Let us describe now the maximal ideal space of the algebra $APH(\Sigma)$.
  
Let $R=\{z \in \mathbb C: e^{-2\pi^2} \leq |z| \leq 1\}$. Then $\pi(z):=e^{2\pi i z}$, $z\in\mathbb C$, determines a projection of the strip $\Sigma$ onto the annulus $R$, so that
the triple $(\Sigma,R,\pi)$ forms a principal bundle on $R$ with fibre $\mathbb Z$. Suppose that $U_1$ and $U_2$ are compact simply connected subsets of $R$ which cover $R$.  Then we can represent $\Sigma$ as a quotient space of $(U_1 \times \mathbb Z) \sqcup (U_2 \times \mathbb Z)$ by
the equivalence relation $\sim$ determined by a locally constant function $c_{12}:U_1 \cap U_2 \mapsto \mathbb Z$ such that
$U_1 \times \mathbb Z \ni (z,n) \sim \bigl(z,n+c_{12}(z)\bigr) \in U_2 \times \mathbb Z$
for all $z \in U_1 \cap U_2$, $n \in \mathbb Z$. 

Let $b\mathbb Z$ be the Bohr compactification of $\mathbb Z$. Then $b\mathbb Z$ is a compact abelian group with $\mathbb Z \subset b\mathbb Z$ acting continuously on $b\mathbb Z$. In particular, given $z \in R$ we have  a continuous mapping $b\mathbb Z \mapsto b\mathbb Z$ determined by the formula
$\xi \mapsto \xi+c_{12}(z)$, $\xi \in b\mathbb Z$.
So, we can define
\begin{equation*}
b\Sigma:=(U_1 \times b\mathbb Z) \sqcup (U_2 \times b\mathbb Z)/\sim
\end{equation*}
where, by definition,
$U_1 \times b\mathbb Z \ni (z,\xi) \sim \bigl(z,\xi+c_{12}(z)\bigr) \in U_2 \times b\mathbb Z$
for all $z \in U_1 \cap U_2$, $\xi \in b\mathbb Z$. The local embeddings $U_k \times \mathbb Z \hookrightarrow U_k \times b\mathbb Z$ ($k=1,2$) determine an embedding $\iota_0:\Sigma \hookrightarrow b\Sigma$. Similarly, an embedding $\iota_\xi:\Sigma \hookrightarrow b\Sigma$, where $\xi \in b\mathbb Z$, is determined by the local embedding $(z,n) \mapsto (z,n+\xi)$. Since $\mathbb Z$ is dense in $b\mathbb Z$, $\iota_\xi(\Sigma)$ is dense in $b\Sigma$ for each $\xi \in b\Sigma$. Furthermore, the sets $\iota_\xi(\Sigma)$ are mutually disjoint and cover $b\Sigma$. 

We say that a function $f:b\Sigma \mapsto \mathbb C$ is \textit{holomorphic} if its pullback $\iota_\xi^* f$ is holomorphic on $\Sigma$ for each $\xi \in b\mathbb Z$. 
The space of functions holomorphic on $b\Sigma$ is denoted by $\mathcal O(b\Sigma)$. Now, we have
 
\begin{proposition}[\cite{BK}]
\label{prop2}
$APH(\Sigma)$ is naturally isomorphic to $\mathcal O(b\Sigma).$ 
\end{proposition}

The result states that $\mathcal O(b\Sigma)|_{\Sigma}=APH(\Sigma)$. Moreover, as it is shown in \cite{BK}, if $\iota_\xi^* f$ is holomorphic on $\Sigma$ for a certain $\xi \in b\Sigma$, then $f$ is holomorphic on $b\Sigma$. This gives another definition of a holomorphic almost periodic function (cf. Definition \ref{def2}). 

\begin{proposition}[\cite{BK}]
\label{prop4}
$b\Sigma$ is homeomorphic to the maximal ideal space of the algebra $APH(\Sigma)$.
\end{proposition}

\noindent {\bf 2.2.} In the proof of Theorem \ref{thm3} we use the following equivalent definition of $B$-valued holomorphic functions on $\Sigma$ (cf. Definition \ref{def2}).

\begin{proposition}\label{exten}
A continuous function $f:\Sigma \mapsto B$ is  holomorphic almost periodic on $\Sigma$ if and only if it is holomorphic in the interior of $\Sigma$ and admits a continuous extension to $b\Sigma$.
\end{proposition}
\begin{proof}
Assume that $f$ is a $B$-valued continuous almost periodic function on $\Sigma$. Since $g\circ f\in APH(\Sigma)$ for every continuous functional $g\in B^*$, the function $f$ admits a continuous extension
$\widehat f:b\Sigma\to B^{**}$; here $B^{**}$ is equipped with weak$^*$ topology. Next, due to Definition \ref{def2},
the closure $\overline{f(\Sigma)}$ of $f(\Sigma)$ in the strong topology of $B$ is compact. We naturally identify $B$ with a subspace of $B^{**}$. Since any compact subset of $B$ is also compact in the weak$^*$ topology of $B^{**}$, $\widehat f(b\Sigma)=\overline{f(\Sigma)}\subset B$. Hence $\widehat f$ maps $b\Sigma$ into $B$ and is continuous.

Conversely, assume that $f:b\Sigma\to B$ is continuous and $f|_{\Sigma}$ is holomorphic. Since $f$ is uniformly continuous on $b\Sigma$ and the natural action of group $\mathbb R$ on $\Sigma$ by shifts is extended to a continuous action of $b\mathbb R$ on $b\Sigma$, the family of shifts
$\{S_x f\}_{x\in\mathbb R}$ is relatively compact in $C_b^B(b\Sigma)\subset C_b^B(\Sigma)$.
\end{proof}

Now, suppose that $f \in APH_B(\Sigma)$. According to Proposition \ref{exten}, there exists a continuous extension $\hat{f}$ of $f$ to $b\Sigma$. For each continuous functional $g\in B^*$ the function $g\circ f$ belongs to $APH(\Sigma)$. By Proposition \ref{prop2}, $g\circ f$ admits a continuous extension $\widehat{g\circ f}\in {\mathcal O}(b\Sigma)$. Since $\Sigma$ is dense in $b\Sigma$, the identity $g\circ\hat f=\widehat{g\circ f}$ is valid for all $g\in B^*$. This implies that for each $\xi\in b\mathbb Z$, the function 
$g\circ\hat f\circ\iota_{\xi}\in APH(\Sigma)$ for all $g\in B^*$. Therefore the continuous $B$-valued function $\hat f\circ\iota_{\xi}$ is holomorphic in the interior of $\Sigma$, i.e., it belongs to $APH_B(\Sigma)$.

\sect{Proof of Theorem 1.4}
\begin{proposition}\label{prop3.1}
$APH(\Sigma)$ has metric approximation property.
\end{proposition}
\begin{proof}
We refer to the book of Besicovich \cite{Bes} for the corresponding definitions and facts from the theory of almost periodic functions.

Let $K \subset APH(\Sigma)$ be compact. Given $\varepsilon>0$ consider an $\frac{\varepsilon}{3}$-net $\{f_1,\dots, f_l\}\subset K$. Let 
\begin{equation*}
\mathcal K(t):=\sum_{|\nu_1| \leq n_1, \dots, |\nu_r| \leq n_r} \left(1-\frac{\nu_1}{n_1} \right)\dots \left(1-\frac{\nu_r}{n_r} \right)e^{-i \left(\frac{\nu_1}{n_1}\beta_1+\dots+\frac{\nu_r}{n_r}\beta_r \right)t}
\end{equation*}
be a Bochner-Fejer kernel such that for all $1\le k\le l$
\begin{equation}\label{eq3.1}
\sup_{z\in\Sigma}|f_k(z)-M_t\{f_k(z+t)\mathcal K(t)\}|\le\frac{\varepsilon}{3}.
\end{equation}
Here $\beta_1,\dots,\beta_r$ are elements of a basis over $\mathbb Q$ of the union of spectra of functions $f_1,\dots, f_l$,\, $\nu_1,\dots, \nu_r\in\mathbb Z$, $n_1,\dots, n_r\in\mathbb N$ and 
$$
M_t\{f_k(z+t)\mathcal K(t)\}:=\lim_{T\to\infty}\frac{1}{2T}\int_{-T}^{T}f_k(z+t)\mathcal K(t)\,dt
$$
are the corresponding holomorphic Bochner-Fejer polynomials.

We define an operator $T: APH(\Sigma)\to APH(\Sigma)$ from Definition \ref{def2'}  by the formula
\begin{equation}\label{eq3.2}
(Tf)(z):=M_t\{f(z+t)\mathcal K(t)\},\ \ \ f\in APH(\Sigma).
\end{equation}
Then $T$ is a linear projection onto a finite-dimensional space generated by functions $e^{i \left(\frac{\nu_1}{n_1}\beta_1+\dots+\frac{\nu_r}{n_r}\beta_r \right)z}$, $|\nu_1| \leq n_1, \dots, |\nu_r| \leq n_r$. Moreover, since $\mathcal K(t)\ge 0$ for all $t\in\mathbb R$ and
$M_t\{\mathcal K(t)\}=1$, the norm of $T$ is $1$. Finally, given $f\in K$ choose $k$ such that $\|f-f_k\|_{APH(\Sigma)}\le\frac{\varepsilon}{3}$. Then we have by \eqref{eq3.1}
\begin{equation*}
\|Tf-f\|_{APH(\Sigma)} \leq \|T(f-f_k)\|_{APH(\Sigma)}+\|Tf_k-f_k\|_{APH(\Sigma)}+\|f_k-f\|_{APH(\Sigma)}<\varepsilon .
\end{equation*}

This completes the proof of the proposition.
\end{proof}

Let us prove now Theorem \ref{thm3}. According to Section 2.2 each $B$-valued almost periodic function $f$ on $\Sigma$ admits a continuous extension $\hat f\in C_b^B(b\Sigma)$. Also, for each $g\in B^*$ the function $g\circ \hat f\in\mathcal O(b\Sigma)$. Since $\mathcal O(b\Sigma)\cong APH(\Sigma)$ has the approximation property by Proposition \ref{prop3.1} and since $b\Sigma$ is compact, the results of Section 5.1 of \cite{Gro} imply that each $\hat f$ belongs to the closure in $C_b^B(b\Sigma)$ of the symmetric tensor product of $APH(\Sigma)$ and $B$. Finally, the classical Bohr theorem, see, e.g., \cite{Bes}, asserts that each element of $APH(\Sigma)$ is the uniform limit of finite linear combinations of functions $e^{i\lambda z}$, $\lambda\in\mathbb R$. These facts complete the proof of the theorem.\ \ \ \ $\Box$

\sect{Proofs of Theorems 1.9 and 1.11}
\subsect{Semi-almost periodic holomorphic functions on $\mathbb D$ with values in a Banach space}
To prove Theorem \ref{approx2} we generalize results of \cite{BK} to the Banach space of semi-almost periodic holomorphic functions on the unit disk $\mathbb D$ with values in a Banach space $B$. 
\\
{\bf 4.1.1.} Let $L_B^\infty(\partial\mathbb D)$ be the Banach space of $B$-valued bounded measurable functions on $\partial\mathbb D$ equipped with $\sup$-norm.

\begin{definition}\label{d5.1}
A function $f \in L_B^\infty(\partial\mathbb D)$ is called semi-almost periodic on $\partial\mathbb D$ if for any $t_0 \in [0,2\pi)$ and any $\varepsilon>0$ there exist a number $s=s(t_0,\varepsilon) \in (0,\pi)$ and functions $f_k:\gamma_{t_0}^k(s)\to B$, $\gamma_{t_0}^k(s):=\{e^{i(t_0+kt)} :\ 0 \leq t<s<2\pi\}$, $k\in\{-1,1\}$, such that functions
\begin{equation*}
t \mapsto f_k\bigl(e^{i(t_0+kse^t)}\bigr), \quad -\infty<t<0, \quad k \in \{-1,1\},
\end{equation*}
are restrictions of some $B$-valued almost periodic functions from $AP_B(\mathbb R)$
and
\begin{equation*}
\sup_{z \in \gamma_{t_0}^k(s)}\|f(z)-f_k(z)\|_B<\varepsilon,\ \ \ k\in\{-1,1\}.
\end{equation*}
\end{definition}

Analogously, we denote by $SAP_B(S) \subset L_B^\infty(S)$ the Banach space of semi-almost periodic functions that are continuous on $\partial\mathbb D \setminus S$. 

Using the Poisson integral formula we can extend each function from $SAP_B(S)$ to a bounded $B$-valued harmonic function on $\mathbb D$ having the same $\sup$-norm. Below we identify $SAP_B(S)$ with its harmonic extension.

Suppose that $S$ contains at least two points. Recall that $A_S$ is the closure in $H^\infty(\mathbb D)$ of the algebra generated by $A(\mathbb D)$ and holomorphic functions of the form $e^{f}$, where $\Real(f)|_{\partial\mathbb D}$ is a finite linear combination (over $\mathbb R$) of characteristic functions of closed arcs in $\partial\mathbb D$ whose endpoints belong to $S$. If $S$ consists of a single point, then $A_S$ is the closure of the algebra generated by $A(\mathbb D)$ and functions $ge^{\lambda f}$, where $\Real(f)|_{\partial\mathbb D}$ is the characteristic function of a closed arc with an endpoint in $S$ and $g \in A(\mathbb D)$ is a function such that $ge^{f}$ has discontinuity on $S$ only. (Here $A({\mathbb D})\subset H^\infty(\mathbb D)$ is the disk-algebra, see Section 1.)

We define the Banach space  $A_S^B$ as 
\begin{equation*}
A_S^B=A_S \otimes B.
\end{equation*}
Here $\otimes$ stands for completion of symmetric tensor product with respect to norm
\begin{equation}
\label{tensnorm}
\left\|\sum_{k=1}^m f_k b_k\right\|:=\sup_{z \in \partial D} \left\|\sum_{k=1}^m f_k(z) b_k\right\|_B\ \ \ {\rm with}\ \ \ f_k\in A_S,\ b_k\in B.
\end{equation}

Let $H^\infty_B(\mathbb D)$ be the Banach space of bounded $B$-valued holomorphic functions on $\mathbb D$  equipped with $\sup$-norm. The Banach-valued analogue of Theorems \ref{approx2} is now formulated as follows.

\begin{theorem}
\label{Bapprox}
$SAP_B(S) \cap H^\infty_B(\mathbb D)=A_S^B$.
\end{theorem}
\noindent {\bf 4.1.2.}
For the proof we require some auxiliary results.

Let $APC(\Sigma)$ be the Banach algebra of functions $f:\Sigma \mapsto \mathbb C$ uniformly continuous on $\Sigma:=\{z\in\mathbb C\ :\ \Imag(z)\in [0,\pi]\}$ and almost periodic on each horizontal line. 

\begin{proposition}[\cite{BK}]
\label{propcoin}
$M(APH(\Sigma))=M(APC(\Sigma))$.
\end{proposition}

We define $APC_B(\Sigma)=APC(\Sigma) \otimes B$.

\begin{lemma}
\label{fivelem}
Suppose that $f_1 \in AP_B(\mathbb R)$, $f_2 \in AP_B(\mathbb R+i\pi)$. Then there exists a function $F \in APC_B(\Sigma)$ which is harmonic in the interior of $\Sigma$ whose boundary values are $f_1$ and $f_2$. Moreover, $F$ admits a continuous extension to the maximal ideal space $M(APH(\Sigma))$.
\end{lemma}
\begin{proof}
By definition $f_1$ and $f_2$ can be approximated on $\mathbb R$ and $\mathbb R +i\pi$ by functions of the form
\begin{equation*}
q_1=\sum_{l=0}^{k_1} b_l e^{i\lambda_l t}, \quad q_2=\sum_{l=0}^{k_2} c_l e^{i\mu_l t},
\end{equation*}
respectively,
where $b_l$, $c_l \in B$, $\lambda_l$, $\mu_l \in \mathbb R$. 
Using the Poisson integral formula we extend $q_1$ and $q_2$ from the boundary to a function $H\in APC_B(\Sigma)$ harmonic in the interior of $\Sigma$, see, e.g., \cite{BK}, Lemma 4.3 for similar arguments. Also, by the definition of the algebra $APC_B(\Sigma)$ and by Proposition \ref{propcoin}, $H$ admits a continuous extension to $M(APH(\Sigma))$. Now, by the 
maximum principle for harmonic functions the sequence of functions $H$ converges in $APC_B(\Sigma)$ to a certain function $F$ harmonic in the interior of $\Sigma$ which satisfies the required properties of the lemma.
\end{proof}

Let $\phi_{z_{0}}:\mathbb D\to\mathbb H_{+}$\,,
\begin{equation}\label{fi}
\phi_{z_{0}}(z):=\frac{2i(z_{0}-z)}{z_{0}+z},\ \ \ z\in\mathbb D,
\end{equation}
be a conformal map of $\mathbb D$ onto the upper half-plane $\mathbb H_{+}$. Then
$\phi_{z_{0}}$ is also continuous on $\partial\mathbb D\setminus\{-z_{0}\}$ and maps it diffeomorphically onto $\mathbb R$ (the boundary of $\mathbb H_{+}$) so that $\phi_{z_{0}}(z_{0})=0$.
Let $\Sigma_{0}$ be the interior of the strip $\Sigma$. Consider the conformal map ${\rm Log}:\mathbb H_{+}\to\Sigma_{0}$,
$z\mapsto {\rm Log}(z):=\ln|z|+i{\rm Arg}(z)$, where ${\rm Arg}:\mathbb C\setminus\mathbb R_{-}\to (-\pi,\pi)$ is the principal branch of the multi-function ${\rm arg}$, the argument of a complex number. The function ${\rm Log}$ is extended to a homeomorphism of $\overline{\mathbb H}_{+}\setminus\{0\}$ onto $\Sigma$; here $\overline{\mathbb H}_{+}$ stands for the closure of $\mathbb H_{+}$.

The proof of the next statement uses Lemma \ref{fivelem} and is very similar to the proof of Lemma 4.2 (for $B=\mathbb C$) in \cite{BK}, so we omit it.

Suppose that $z_{0}=e^{it_{0}}$. For $s\in (0,\pi)$ we set $\gamma_{1}(z_0,s):={\rm Log}(\phi_{z_{0}}(\gamma_{t_{0}^{1}}(s)))\subset\mathbb R$ and 
$\gamma_{-1}(z_0,s):={\rm Log}(\phi_{z_{0}}(\gamma_{t_{0}^{-1}}(s)))\subset \mathbb R+i\pi$. 

\begin{lemma}
\label{fourlem}
Let $z_0 \in S$, suppose that $f \in SAP_B(\{-z_0,z_0\})$. We put $f_k=f|_{\gamma_{t_0^k}(\pi)}$, and define on arc $\gamma_k(z_0,s)$ $$h_k=f_k \circ \varphi_{z_0}^{-1} \circ \Log^{-1}, \quad k \in \{-1,1\}.$$ Then for any $\varepsilon>0$ there exist a number $s_\varepsilon \in (0,s)$ and a function $H\in APC_B(\Sigma)$ harmonic on $\Sigma_0$ such that 
\begin{equation*}
\sup_{z \in \gamma_k(z_0,s_\varepsilon)}\|h_k(z)-H(z)\|_B<\varepsilon, \quad k \in \{-1,1\}.
\end{equation*}
\end{lemma}

Let $z_0 \in \partial\mathbb D$ and $U_{z_0}$ be the intersection of an open disk of radius $ \leq 1$
centered at $z_0$ with $\bar{\mathbb D} \setminus z_0$. We call such $U_{z_0}$ a {\em circular neighbourhood} of $z_0$.

We say that a bounded continuous function $f:\mathbb D \mapsto B$ is {\em almost-periodic near} $z_0$ if there exist a circular neighbourhood $U_{z_0}$, and a function $\hat{f} \in APC_B(\Sigma)$ such that
\begin{equation*}
f(z)=\hat{f}\bigl(\Log\bigl(\varphi_{z_0}(z)\bigr)\bigr), \quad z \in U_{z_0}.
\end{equation*}

Let $M_S$ be the maximal ideal space of $H^\infty(\mathbb D)\cap SAP(S)$ and $a_S:M_S\to\bar{\mathbb D}$ be the natural continuous surjection, see Remark \ref{struct}. Recall that for $z_0\in S$, $a_S^{-1}(z_0)$ can be naturally identified with $b\Sigma:=M(APH(\Sigma))$. Therefore the algebra ${\mathcal O}(a_S^{-1}(z_0))$ of holomorphic functions on $a_S^{-1}(z_0)$ can be defined similarly to ${\mathcal O}(b\Sigma)$ from Section 2 (by means of this identification), see \cite{BK} for details. 

In the proof of Theorem 1.8 of \cite{BK} (see Lemmas 4.4, 4.6 there) we established
\begin{itemize}
\item[(1)]
Any scalar harmonic function $f$ on $\mathbb D$ almost periodic near $z_0$ admits a continuous extension to $a_S^{-1}(\bar{U}_{z_0})\subset M_S$ for some circular neighbourhood $U_{z_{0}}$.
\item[(2)]
For any holomorphic function $f\in {\mathcal O}(a_S^{-1}(z_0))$,
$z_0\in S$, there is a bounded holomorphic function $\hat f$ on $\mathbb D$ of the same $\sup$-norm almost periodic near $z_0$ such that its extension to $a_S^{-1}(z_0)$ coincides with $f$.
\end{itemize}

Similarly to ${\mathcal O}(b\Sigma)$ we define the Banach space 
${\mathcal O}_{B}(b\Sigma)$ of $B$-valued holomorphic functions on $b\Sigma$, cf. Section 2. Then by Theorem \ref{thm3}, ${\mathcal O}_B(b\Sigma):={\mathcal O}(b\Sigma)\otimes B$. This together with the $B$-valued Bohr approximation theorem for the subalgebra of harmonic functions in $APC_B(\Sigma)$ (see Lemma \ref{fivelem}) and statements (1), (2) imply
\begin{itemize}
\item[(3)]
Any $B$-valued harmonic function $f$ on $\mathbb D$ almost periodic near $z_0$ admits a continuous extension to $a_S^{-1}(\bar{U}_{z_0})\subset M_S$ for some circular neighbourhood $U_{z_{0}}$.
\item[(4)]
For any holomorphic function $f\in {\mathcal O}_B(a_S^{-1}(z_0))$,
$z_0\in S$, there is a bounded $B$-valued holomorphic function $\hat f$ on $\mathbb D$ of the same $\sup$-norm almost periodic near $z_0$ such that its extension to $a_S^{-1}(z_0)$ coincides with $f$.
\end{itemize}

From (3) and (4) we obtain
\begin{lemma}\label{l5.6}
Let $f\in SAP_B(S)\cap H_B^{\infty}(\mathbb D)$ and $z_0\in\partial\mathbb D$. There is a bounded $B$-valued holomorphic function $\hat f$ on $\mathbb D$ almost periodic near $z_0$ such that for any $\varepsilon>0$
there is a circular neighbourhood $U_{z_0;\varepsilon}$ of $z_0$ so that
$$
\sup_{z\in U_{z_0;\varepsilon}}||f(z)-\hat f(z)||_B<\varepsilon .
$$
\end{lemma}
\begin{proof}
Assume, first, that $z_0\in S$. By Lemma \ref{fourlem}, for any $n\in\mathbb N$ there exist a number $s_n \in (0,s)$ and a function $H_n\in APC_B(\Sigma)$ harmonic on $\Sigma_0$ such that 
\begin{equation}\label{e5.3}
\sup_{z \in \gamma_k(z_0,s_n)}\|f_k(z)-H_n(z)\|_B<\frac{1}{n}\, , \quad k \in \{-1,1\}.
\end{equation}
Using the Poisson integral formula for the bounded $B$-valued harmonic function $f-H_n$ on $\mathbb D$ we easily obtain from \eqref{e5.3} that there is a circular neighbourhood $V_{z_0;n}$ of $z_0$ such that
\begin{equation}\label{e5.4}
\sup_{z \in V_{z_0;n}}\|f(z)-H_n(z)\|_B<\frac{2}{n} .
\end{equation}
According to (3) each $H_n$ admits a continuous extension $\hat H_n$ to $a_S^{-1}(z_0)\cong b\Sigma$. Moreover, \eqref{e5.4} implies that the restriction of the sequence $\{\hat H_n\}_{n\in\mathbb N}$ to $a_S^{-1}(z_0)$ forms a Cauchy sequence in $C_B(a_S^{-1}(z_0))$. Let $\hat H\in C_B(a_S^{-1}(z_0))$ be the limit of this sequence. 

Further, for any functional $\phi\in B^*$ the function $\phi\circ f\in SAP(S)\cap H^\infty(\mathbb D)$ and therefore admits a continuous extension $f_{\phi}$ to $a_S^{-1}(z_0)$ such that on $a_S^{-1}(z_0)$ the extended function belongs to 
${\mathcal O}(a_S^{-1}(z_0))$. Now, \eqref{e5.4} implies directly that $f_\phi=\phi\circ\hat H$ for any $\phi\in B^*$. Then from the definition of ${\mathcal O}_B(a_S^{-1}(z_0))$, see Section 2, follows that $\hat H\in{\mathcal O}_B(a_S^{-1}(z_0))$. Thus by (4) we find a bounded $B$-valued holomorphic function $\hat f$ on $\mathbb D$ of the same $\sup$-norm as $\hat H$ almost periodic near $z_0$ such that its extension to $a_S^{-1}(z_0)$ coincides with $\hat H$. Now by the definition of the topology of $M_S$, see \cite{BK}, Lemma 4.4 (a), we obtain that for any $\varepsilon>0$ there is a number $N\in\mathbb N$ such that for all $n\geq N$,
$$
\sup_{z \in V_{z_0;n}}\|\hat f(z)-H_n(z)\|_B<\frac{\varepsilon}{2}.
$$ 
Finally, choose $n\geq N$ in \eqref{e5.4} such that the right-hand there $<\frac{\varepsilon}{2}$. For this $n$ we set $U_{z_0;\varepsilon}:=V_{z_0;n}$. Then the previous inequality and \eqref{e5.4} imply the required
$$
\sup_{z\in U_{z_0;\varepsilon}}||f(z)-\hat f(z)||_B<\varepsilon .
$$

If, now, $z_0\not\in S$, then, by definition, $f|_{\partial\mathbb D}$ is continuous at $z_0$. In this case as the function $\hat f$ we can choose the constant $B$-valued function equal to $f(z_0)$ on $\mathbb D$. Then the required result follows from the Poisson integral formula for $f-\hat f$.
We leave the details to the reader.
\end{proof}
\noindent {\bf 4.1.3.} We are now ready to prove Theorem \ref{Bapprox}.
\begin{proof}
As it was shown in \cite{BK}, $A_S \subset SAP(S)\cap H^{\infty}(\mathbb D)$. So
$A^B_S \subset {SAP}_B(S) \cap H_B^\infty(\mathbb D)$.
Let us prove the opposite inclusion.

(A) Consider first the case $S=F$, where $F=\{z_i\}_{i=1}^m$ is a finite subset of $\partial\mathbb D$. 

Let $f \in SAP_B(F) \cap H^\infty_B(\mathbb D)$. Then according to Lemma \ref{l5.6} there exists a function $f_{z_1} \in APH_B(\Sigma)$
such that the bounded $B$-valued holomorphic function $g_{z_1}-f$, $g_{z_1}:=f_{z_1}\circ\Log\circ\varphi_{z_0}$, on $\mathbb D$ is continuous and equals $0$ at $z_1$.

Let us show that $g_{z_1} \in A^B_{\{z_1,-z_1\}}$. Indeed, since $f_{z_1} \in APH_B(\Sigma)$, by Theorem \ref{thm3} it can be approximated in $APH_B(\Sigma)$ by finite sums of functions $b e^{i\lambda z}$, $b \in B$, $\lambda \in \mathbb R$, $z \in \Sigma$. In turn, $g_{z_1}$ can be approximated by finite sums of functions $b e^{i\lambda \Log \circ \varphi_{z_1}}$. As it was shown in \cite{BK}, $e^{i\lambda \Log \circ \varphi_{z_1}} \in A_{\{z_1,-z_1\}}$. Hence, $g_{z_1} \in A^B_{\{z_1,-z_1\}}$.

We define
\begin{equation*}
\hat{g}_{z_1}=\frac{g_{z_1}(z)(z+z_1)}{2z_1}.
\end{equation*}
Then, since the function $z \mapsto (z+z_1)/(2z_1) \in A(\mathbb D)$ and equals $0$ at $-z_1$, and $g_{z_1} \in A^B_{\{z_1,-z_1\}}$, the function $\hat{g}_{z_1} \in A^B_{\{z_1\}}$. Moreover, by the definitions of $g_{z_1}$ and $\hat{g}_{z_1}$, the function $\hat{g}_{z_1}-f$ is continuous and equal to zero at $z_1$.  Thus, 
$$
\hat{g}_{z_1}-f \in SAP(F \setminus \{z_1\})\cap H_B^\infty(\mathbb D) .
$$ 
We proceed in this way to get functions $\hat{g}_{z_k} \in A_{\{z_k\}}^B$ such that 
\begin{equation*}
f -\sum_{k=1}^m \hat{g}_{z_k} \in A_B(\mathbb D).
\end{equation*}
Here $A_B(\mathbb D)$ is the Banach space of $B$-valued bounded holomorphic functions on $\mathbb D$ continuous up to the boundary. As in the scalar case using the Taylor expansion at $0$ of functions from $A_B(\mathbb D)$ one can easily show that $A_B(\mathbb D)=A(\mathbb D)\otimes B$. This and the above implication imply that $f\in A_F^B:=A_F\otimes B$, as required.

(B) Let us consider the general case of $S \subset \partial\mathbb D$ an arbitrary closed set. 

Let $f \in SAP_B(S) \cap H^\infty_B(\mathbb D)$. As it follows from Lemma \ref{l5.6} and the arguments presented in part (A), given an $\varepsilon>0$ there exist points $z_k \in \partial\mathbb D$, functions $f_k \in A_{\{z_k\}}^B$ and circular neighbourhoods $U_{z_k}$ ($1 \leq k \leq m$) such that $\{U_k\}_{k=1}^m$ forms an open cover of set $\partial\mathbb D \setminus \{z_k\}_{k=1}^m$ and
\begin{equation}
\label{b3}
\|f(z)-f_k(z)\|_B<\varepsilon\ \ \text{ on }\ \ U_{z_k}, \quad 1 \leq k \leq m.
\end{equation}
Since $S$ is closed, for $z_k \not\in S$ we may assume that $f_k$ is continuous in $\bar{U}_{z_k}$.

Let us define a $B$-valued $1$-cocycle $\{c_{kj}\}_{k,j=1}^m$ on intersections of the sets in $\{U_{z_k}\}_{k=1}^m$ by the formula
\begin{equation}\label{e5.6'}
c_{kj}(z):=f_k(z)-f_j(z), \quad z \in U_{z_k} \cap U_{z_j}.
\end{equation}
Then \eqref{b3} implies
\begin{equation}\label{e5.6}
\sup_{k,j,z}||c_{kj}(z)||_B<2\varepsilon .
\end{equation}

Since $m<\infty$, we may assume, without loss of generality, that none of the intersections $U_{z_k} \cap U_{z_j}$, $k \ne j$, contains points $z_k$. Further, we may choose the above functions and sets such that $c_{kj}$ is holomorphic in $U_{z_k} \cap U_{z_j}$ and continuous in the closure of $U_{z_k} \cap U_{z_j}$. Let $\{\rho\}_{k=1}^m$ be a smooth partition of unity subordinate to the open cover $\{U_{z_k}\}_{k=1}^m$ of an open annulus $A\subset\subset\cup_{k=1}^m U_{z_k}$ with outer boundary $\partial\mathbb D$ such that each $\rho_k$ is the restriction of a smooth function defined on a neighbourhood of $\partial{\mathbb D}$ and $\rho_k(z_k)=1$, $1\leqslant k \leqslant m$. We resolve the cocycle $\{c_{kj}\}_{k,j=0}^m$ using this partition of unity by formulas
\begin{equation*}
\tilde{f}_j(z)=\sum_{k=1}^m \rho_k(z)c_{kj}(z), \quad z \in U_{z_j}\cap A, 
\end{equation*}
so, by definition,
\begin{equation}
\label{cij}
c_{kj}(z)=\tilde{f}_k(z)-\tilde{f}_j(z), \quad z \in U_{z_k} \cap U_{z_j}\cap A.
\end{equation}
Now since $c_{kj}$ are $B$-valued holomorphic functions in $U_{z_i} \cap U_{z_j}$ continuous up to the boundary, $\rho_{kj}$ are smooth functions on a neighbourhood of $\partial{\mathbb D}$ and $A\subset\subset\cup_{k=1}^m U_{z_k}$,
$$
h(z):=\frac{\partial \tilde{f}_j(z) }{\partial \bar{z}},\ \ \ z\in U_{z_{j}}\cap A,\ \ \ 1\leq j\leq m,
$$
is a smooth global $B$-valued function on $A$ continuous up to the boundary.

We define
\begin{equation}
\label{bint}
H(z):=\frac{1}{2\pi i} \int_{\zeta \in \mathbb D} \frac{h(\zeta)}{\zeta-z} d\zeta \wedge d\bar{\zeta}\, ,\ \ \
z \in \bar{\mathbb D}. 
\end{equation}
Similarly to the scalar case we use polar coordinates in (\ref{bint}) to obtain
\begin{equation*}
\sup_{z \in A}\|H(z)\|_B \leqslant C w(A)\sup_{z \in A} \|h(z)\|_B,
\end{equation*}
where $w(A)$ is width of $A$ and $C>0$ is a numerical constant.
Furthermore, $H$ is continuous on $\bar{A}$ and $\frac{\partial H}{\partial \bar{z}}=h$ on $A$. 
Without loss of generality (see \eqref{e5.6} and the definition of functions $\rho_j$) we may assume that $w(A)$ is sufficiently small, so that
\begin{equation*}
\sup_{z \in A} \|H(z)\|_B<\varepsilon.
\end{equation*}

Let us define
\begin{equation}
\label{cfh}
c_i(z):=\tilde{f}_i(z)-H(z), \quad z \in \bar{U}_{z_i}\cap\bar{A}.
\end{equation}
Then $c_i$ is a $B$-valued continuous function holomorphic in $U_{z_i}\cap A$ satisfying (see \eqref{e5.6})
\begin{equation}\label{e5.12'}
\sup_{z\in U_{z_i}\cap A}||c_i(z)||_B\leq 3\varepsilon.
\end{equation}
By (\ref{cij}), 
\begin{equation}\label{e5.12}
c_i(z)-c_j(z)=c_{ij}(z), \quad z \in U_{z_i} \cap U_{z_j}\cap A.
\end{equation}

Let us determine a function $f_\varepsilon$ defined on $\bar{A} \setminus \{z_i\}_{i=1}^m$ by formulas
\begin{equation*}
f_\varepsilon(z):=f_i(z)-c_i(z), \quad z \in U_{z_i}\cap\bar{A}.
\end{equation*}
According to \eqref{e5.6'} and \eqref{e5.12}, $f_{\varepsilon}$ is a bounded continuous $B$-valued function on $\bar{A}\setminus \{z_i\}_{i=1}^m$ holomorphic in $A$.
Furthermore, since $c_i$ is continuous on $\bar{U}_{z_i}\cap\bar{A}$, and $f_i\in A_{\{z_i\}}^B$, for $z_i\in S$ and $f_i\in A_B(\mathbb D)$ otherwise,  $f_\varepsilon|_{\partial\mathbb D}\in SAP_B(F)$ where $F=\{z_i\}_{i=1}^m \cap S$. Also, from inequalities (\ref{b3}) and \eqref{e5.12'} we have
\begin{equation}
\label{ineq1t}
\sup_{z\in A}\|f(z)-f_\varepsilon(z)\|_B<4\varepsilon .
\end{equation}

Now let $D'$ be an open disk centered at $0$ whose intersection with $A$ is an annulus of width $\frac{w(A)}{2}$. Let $A'$ be the open annulus with the interior boundary coinciding with the interior boundary of $A$ and with the outer boundary $\{z\in\mathbb C\ :\ |z|=2\}$. Then $\{A',D'\}$ forms an open cover of $\mathbb D$. Let us define a $1$-cocycle on $A' \cap D'=A \cap D'$ by the formula
\begin{equation*}
c(z)=f(z)-f_\varepsilon(z), \quad z \in A \cap D'.
\end{equation*}
As it follows from (\ref{ineq1t}), $\|c(z)\|_B<4\varepsilon$ for all $z \in D' \cap A'$.
Consider a smooth partition of unity subordinate to the cover $\{A',D'\}$ of $\mathbb D$ which consists of smooth radial functions $\rho_{A'}$, $\rho_{D'}$ defined on $\mathbb C$ such that
\begin{equation}
\label{coolineq}
\max\{\|\nabla \rho_{A'}\|_{L^\infty(\mathbb C)}, \|\nabla \rho_{D'}\|_{L^\infty(\mathbb C)}\} \leqslant \frac{C}{w(D' \cap A')}=\frac{2C}{w(A)}
\end{equation}
where $C>0$ is a numerical constant. We resolve the cocycle $c$ as follows:
$$
\begin{array}{l}
f_{A'}(z)=-\rho_{D'}(z)c(z), \quad z \in A',\\
f_{D'}(z)=\rho_{A'}(z)c(z), \qquad z \in D'.
\end{array}
$$
Then $c(z)=f_{D'}(z)-f_{A'}(z)$, $z \in A' \cap D'$. Since $c$ is holomorphic in $A' \cap D'$, 
\begin{equation}
g(z)=\left\{
\begin{array}{cc}
\displaystyle
\frac{\partial f_{A'}(z)}{\partial \bar{z}},&z \in A',
\\
[4mm]
\displaystyle\frac{\partial f_{D'}(z)}{\partial \bar{z}},&z \in D'
\end{array}
\right.
\end{equation}
is a bounded smooth $B$-valued function on $\mathbb D$ with support in $A$. 

Next, 
\begin{equation}
\label{d3}
G(z)=\frac{1}{2\pi i} \int_{\zeta \in D} \frac{g(\zeta)}{\zeta-z} d\zeta \wedge d\bar{\zeta}
\end{equation}
is a smooth $B$-valued function on $\mathbb D$, continuous up to the boundary such that $
\frac{\partial G}{\partial \bar{z}}=g$ on $\mathbb D$. Moreover,
by \eqref{ineq1t}, \eqref{coolineq} and the fact that ${\rm supp}(g)\subset A$ we obtain from \eqref{d3}
\begin{equation}\label{e5.17}
\sup_{z\in\mathbb D}\|G(z)\|_B\leq C' w(A)\frac{\varepsilon}{w(A)}=C'\varepsilon
\end{equation}
for a numerical constant $C'>0$.

Now we define
\begin{equation*}
c_{A'}(z)=f_{A'}(z)-G(z), \quad z \in A,
\end{equation*}
\begin{equation*}
c_{D'}(z)=f_{D'}(z)-G(z), \quad z \in D'.
\end{equation*}
Then $c_{A'}$ and $c_{D'}$ are $B$-valued holomorphic functions in $A$ and $D'$, respectively. Clearly, we have $c_{D'}(z)-c_{A'}(z)=c(z)$ for all $z \in D' \cap A'$. Furthermore, as it follows from (\ref{ineq1t}), \eqref{e5.17} there exists a numerical constant $\bar{C}>0$ such that
\begin{equation*}
\|c_{A'}(z)\|_B<\bar{C}\varepsilon, \quad z \in A,
\quad \|c_{D'}(z)\|_B<\bar{C}\varepsilon, \quad z \in D'.
\end{equation*}

Finally, let us define
\begin{equation*}
F_\varepsilon(z):=\left\{
\begin{array}{ll}
f(z)-c_{D'}(z), \quad z \in D', \\ [1mm]
f_\varepsilon(z)-c_{A'}(z), \quad z \in A. 
\end{array}
\right.\\ [2mm]
\end{equation*}
Clearly, $F_\varepsilon$ is a $B$-valued holomorphic function on $\mathbb D$ (since for every $z \in D' \cap A$ we have $f(z)-f_\varepsilon(z)-c_{D'}(z)+c_{A'}(z)=f(z)-f_\varepsilon(z)-c(z)=0$). Also, 
\begin{equation*}
\sup_{z\in\mathbb D}\|f(z)-F_\varepsilon(z)\|_B<\hat{C}\varepsilon 
\end{equation*}
for a numerical constant $\hat{C}>0$, and by definition $F_\varepsilon \in SAP_B(F) \cap H^\infty_B(\mathbb D)$, where $F=\{z_1,\dots,z_m\} \cap S$. 

The last inequality and part (A) of the proof show that the complex vector space generated by spaces $A_F^B$ for all possible finite subsets $F\subset S$ is dense in $SAP_B(S)\cap H_B^\infty(\mathbb D)$. Since the closure of all such $A_F^B$ is $A_S^B$, we
obtain the required: $SAP_B(S)\cap H_B^\infty(\mathbb D)=A_S^B$.
\end{proof}
\begin{remark}\label{r4.7}
From the proof of Theorem \ref{Bapprox} one obtains also that each function from $SAP_B(S)\cap H_B^\infty(\mathbb D)$ admits a continuous extension to the maximal ideal space of the algebra $SAP(S)\cap H^\infty(\mathbb D)$. Then this theorem and the results of Section 5.1 of \cite{Gro} imply that $SAP(S)\cap H^\infty(\mathbb D)$ has the approximation property.
\end{remark}
\subsect{Proofs}
\begin{proof}[Proof of Theorem \ref{approx2}]
Let $S=\prod_{k=1}^nS_k\subset(\partial\mathbb D)^n$ where
$S_k\subset\partial\mathbb D$ is a closed set. We must prove that $H^\infty(\mathbb D^n)\cap SAP^n(S)=A_S^n$.

We prove this statement by induction over $n$. For $n=1$ the identity was already proved. Suppose that it is true for $n-1$, that is, $H^\infty(\mathbb D^{n-1})\cap SAP^{n-1}(S')=A^{n-1}_{S'}$ where $S':=\prod_{k=1}^{n-1} S_k$, and prove it for $n$. Let $S=S' \times S_n$.
\begin{lemma} 
\label{proprepr}
$$SAP^n(S)=SAP_B(S_n)$$ for $B=SAP^{n-1}(S')$ and $$A^n_S=A_{S_n}^C$$ for $C=A^{n-1}_{S'}$.
\end{lemma}
\begin{proof}
The proof follows directly from the definitions of $SAP^n(S)$ and $A^n_S$.
\end{proof}

Let $f\in H^\infty(\mathbb D^n)\cap SAP^n(S)$. According to Lemma \ref{proprepr}, $f\in SAP_B(S_n)$ for $B=SAP^{n-1}(S')$. Also, by the Poisson integral formula for bounded polyharmonic functions on $\mathbb D^n$, and the fact that $f\in H^\infty(\mathbb D^n)$ we obtain that $f\in H_B^{\infty}(\mathbb D)$. Thus, in virtue of Theorem \ref{Bapprox},
\begin{equation}\label{e6.1}
f\in SAP_B(S_n) \cap H^\infty_B(\mathbb D)=A^B_{S_n}.
\end{equation}

Let $\xi=(\xi_1,\dots,\xi_n)$ be coordinates on $(\partial\mathbb D)^n$. Applying again the Poisson integral formula for bounded polyharmonic functions on $\mathbb D^n$   we obtain that the functions $f(\cdot,\dots,\cdot,\xi_n)\in H^{\infty}(\mathbb D^{n-1})$ for almost all $\xi_n\in\partial\mathbb D$. Thus, by the induction hypothesis, for almost all $\xi_n\in\mathbb D$ we have $f(\cdot,\dots,\cdot,\xi_n)\in H^{\infty}(\mathbb D^{n-1})\cap B=A_{S'}^{n-1}$. From here and \eqref{e6.1} we get $f\in A_{S_n}^C$, $C:=A_{S'}^{n-1}$. Now, Lemma \ref{proprepr} implies that $f\in A_S^n$. This shows that $H^\infty(\mathbb D^n)\cap SAP^n(S)\subset A_S^n$. The converse embedding
$A_S^n\subset H^\infty(\mathbb D^n)\cap SAP^n(S)$ follows immediately from the definition of $A_S^n$.
\end{proof}
\begin{proof}[Proof of Theorem \ref{max}]
According to Theorem \ref{approx2} and the scalar version of Theorem \ref{Bapprox}, the algebra $H^\infty(\mathbb D^n)\cap SAP^n(S)$ is symmetric tensor product of algebras $H^{\infty}(\mathbb D)\cap SAP(S_k)$, $1\leq k\leq n$. Therefore by a standard result of the theory of Banach algebras 
we obtain that $M_S$ is homeomorphic to $M_{S_1}\times\cdots\times M_{S_n}$ for the corresponding maximal ideal spaces of these algebras. As it was proved in \cite{BK}, $\mathbb D$ is dense in each $M_{S_k}$ in the corresponding Gelfand topology. Therefore $\mathbb D^n$ is dense in $M_S$ in the Gelfand topology. 
\end{proof}


\begin{thebibliography}{1}

\bibitem{Bes}
A.\,S.\, {Besicovich},
\newblock {\em Almost periodic functions}.
\newblock Dover Publications, 1958.

\bibitem{BK}
A.~{Brudnyi} and D.~{Kinzebulatov},
\newblock On uniform subalgebras of ${L}^\infty$ on the unit circle generated
  by almost periodic functions.
\newblock {Algebra and Analysis} {\bf 19} (2007), 1--33.

\bibitem{E}
P.~{Enflo}, A counterexample to the approximation property in Banach spaces. \newblock {Acta Math.} {\bf 130} (1973), 309--317.

\bibitem{Gar}
J.~{Garnett},
\newblock {\em Bounded analytic functions}.
\newblock Academic Press, 1981.

\bibitem{Gro}
A.~{Grothendieck},
\newblock Products tensoriels toplogiques et espaces nucl\'{e}aires.
\newblock {\em Memoirs Amer. Math. Society} {\bf 16}, 1955.

\bibitem{Sar}
D.~Sarason,
\newblock Toepliz operators with semi-almost periodic kernels.
\newblock {Duke Math J.} {\bf 44} (1977), 357--364.

\end{thebibliography}

\end{document}